\documentclass{elsarticle}

\usepackage[english]{babel}
\usepackage{dsfont} 
\usepackage{graphicx}
\usepackage{color}
\usepackage[hypertexnames=false]{hyperref} 
\usepackage{enumerate}
\usepackage{amsmath}
\usepackage{amssymb,empheq} 
\usepackage{amsthm} 
\usepackage{thmtools} 
\usepackage{cleveref}
\usepackage{xcolor}
\usepackage{subcaption}

\allowdisplaybreaks
\newtheorem{thm}{Theorem}[section]

\newtheorem{lem}[thm]{Lemma}

\newtheorem{cor}[thm]{Corollary}

\theoremstyle{definition}
\newtheorem{definition}[thm]{Definition}

\newtheorem{exm}[thm]{Example}
\newtheorem{exms}[thm]{Examples}

\theoremstyle{remark}
\newtheorem{rem}[thm]{Remark}

\newcommand{\exmsymbol}{\hfill$\circ$}

\newcommand{\cset}{\mathds{C}}
\newcommand{\kset}{\mathds{K}}
\newcommand{\nset}{\mathds{N}}

\newcommand{\rset}{\mathds{R}}

\newcommand{\conv}{\mathrm{conv}\,}

\newcommand{\diff}{\mathrm{d}}

\newcommand{\pos}{\mathrm{Pos}}

\newcommand{\one}{\mathds{1}}

\newcommand{\Bor}{\mathrm{Bor}}
\newcommand{\Herm}{\mathrm{Herm}}
\newcommand{\im}{\mathrm{im}}
\newcommand{\Pos}{\pos}
\newcommand{\tspan}{\mathrm{span}}

\newcommand{\cA}{\mathcal{A}}

\newcommand{\cH}{\mathcal{H}}

\newcommand{\cL}{\mathcal{L}}

\newcommand{\cX}{\mathcal{X}}

\author{
	Lars-Luca Langer\footnote{Email: lars-luca.langer@uni-konstanz.de}$^{\text{,a}}$
}
\address{$^{\text{a}}$Department of Mathematics and Statistics, University of Konstanz, Universit\"atsstra{\ss}e 10, D-78464 Konstanz, Germany}

\journal{ArXiv}

\title{
	A Generalization of Haviland's Theorem for Tensor Spaces
}

\begin{document}

\begin{abstract}
	We generalize Haviland's theorem to abstract polynomials which are given as finite sums of some vector space elements and polynomials. 
	This also generalizes the work of Cimprič and Zalar and their Haviland theorem for matrix or operator polynomials. 
	Unlike Cimprič and Zalar, we do not require a $C^*$-algebra structure but merely require our vector space to be a preordered space to make sense of positivity.
\end{abstract}


\maketitle

\tableofcontents

\bigskip
\hrule height 1.15 pt


\section{Introduction}

	Haviland's theorem (\Cref{cor:HavilandTheorem:Polynomials}) in \cite{havila36} characterizes positive functionals $L: \rset[x_1, \dots, x_n] \to \rset$ on the space of polynomials in $n \in \nset$ variables. 
	The theorem states, that a functional is positive, if and only if, there exists a moment measure $\mu$ on $\rset^n$ such that 
\[
	L(p) = \int_{\rset^n} p(x) ~\diff\mu(x)
\]
	for all $p \in \rset[x_1, \dots, x_n]$. 
	Cimprič and Zalar generalized this Haviland theorem to the case of operator or matrix polynomials  (\Cref{cor:HavilandTheorem:OperatorPolynomials}) in \cite{cimpzalar13}. 
	These operator or matrix polynomials $p \in B(\cH)_h \otimes \rset[x_1, \dots, x_n]$, 
\[
	p = p_1 \otimes A_1 + \dots + p_k \otimes A_k
\]
	are simply finite sums of real polynomials $p_i \in \rset[x_1, \dots, x_n]$ multiplied by operators $A_i \in B(\cH)_h$ where $B(\cH)_h$ is the real Banach space of bounded self-adjoint linear operators on some Hilbert space $\cH$. 
	If the Hilbert space $\cH$ is finite-dimensional, then all these operators are Hermitian matrices and the polynomials in $B(\cH)_h \otimes \rset[x_1, \dots, x_n]$ are matrix polynomials.

	Now, one can question whether a type of Haviland's theorem exists with other positivity cones, not necessarily given by the $C^*$-algebra structure of $B(\cH)_h$. 
	It is easy to see, that Haviland's theorem generalizes to the case where we regard a matrix $M = (m_{i, j})_{i, j = 1, \dots, k}$ as positive if all its components $m_{i, j}$ are positive. 
	Hence, this paper presents a generalization of Haviland's theorem to tensor spaces $V \otimes \rset[x_1, \dots, x_n]$ (\Cref{thm:Haviland:TensorCones}) where $V$ is merely a preordered vector space, dropping the necessity for an algebraic structure. 
	This setting then includes both Cimprič and Zalar's theorem, as well as the case of componentwise positive matrices.

	The paper is structured as follows. 
	First, in \Cref{sec:PreliminaryTheory}, we give some preliminary theory. 
	There, we explain what ordered spaces are and how they relate to convex cones (\Cref{subsec:ConvexConesPositiveFunctionals}), what type of Riesz--Markov--Kakutani theorem we use (\Cref{subsec:RieszMarkovKakutani}), how to integrate vector-valued functions (\Cref{subsec:VectorValuedIntegrationPettisTheorem}), and what adapted spaces are and how they are used to generalize Haviland's theorem (\Cref{subsec:AdaptedSpacesChoquetBasicRepresentationTheorem}). 
	Then, in \Cref{sec:MapValuedMeasures}, we explain what map-valued measures are and what properties they exhibit.
	The definition of map-valued measures is based on the definition of non-negative operator-valued measures from \cite[Sec.\ 2]{cimpzalar13}.
	Lastly, in \Cref{sec:HavilandTheoremTensorSpaces}, we give our main result, a Haviland's theorem for tensor spaces which is based on the Haviland's theorem from \cite[Rhm.\ 3]{cimpzalar13} but removes some redundancies.


\section{Preliminary Theory}
\label{sec:PreliminaryTheory}


\subsection{Convex Cones and Positive Functionals}
\label{subsec:ConvexConesPositiveFunctionals}

\begin{definition}
	Let $V$ be a real vector space. 
	A set $C \subseteq V$ is called a
\begin{itemize}
\item
	\emph{cone}, if $\lambda c \in C$ for all $c \in C$ and $\lambda \in \rset_{\geq 0}$.
\item
	\emph{convex cone}, if $C$ is a cone and $C + C \subseteq C$.
\item
	\emph{salient convex cone}, if $C$ is a convex cone and $C \cap (-C) = \{0\}$.
\item
	\emph{full-dimensional convex cone}, if $C$ is a convex cone and $V = C - C$.
\end{itemize}
\end{definition}

	Every preordered vector space $(V, \geq)$ gives rise to a convex cone 
\[
	V_+ := \{v \in V \mid v \geq 0\}
\]
	called the \emph{positivity cone of $V$}. 
	Conversely, every convex cone $C \subseteq V$ gives rise to a preordered vector space $(V, \geq_C)$ where 
\[
	x \geq_C y \iff x - y \in C.
\]
	Hence, there is a one-to-one correspondence between preordered vector spaces and convex cones. 
	There is also a one-to-one correspondence between ordered vector spaces and salient convex cones, as well as a one-to-one correspondence between directed preordered vector spaces and full-dimensional convex cones. 
	For more on convex cones and preordered vector spaces, see e.g. \cite{schaef99}.

\begin{exms}
\label{exms:OrderedVectorSpaces}
	In the following, we give some examples of preordered vector spaces.
\begin{enumerate}[(a)]
\item
	The space of Hermitian matrices:

	Let $n \in \nset$ and $\Herm_n := \{A \in \cset^{n \times n} \mid \overline{A} = A^T\}$ be the subspace of $n \times n$ Hermitian matrices. 
	Then $\Herm_n$ is a directed ordered vector space with the Loewner order, i.e., 
\[
	A \geq B \iff A - B \text{ is positive semidefinite}.
\]
	That is, the convex cone of positive semidefinite matrices is salient and full-dimensional.
\item
	The space of self-adjoint operators:

	Let $\cH$ be a Hilbert space and $B(\cH)_h$ be the space of linear bounded self-adjoint operators on $\cH$. 
	Similarly to the Hermitian matrices, $B(\cH)_h$ is also a directed ordered vector space with the Loewner order.
	That is, the convex cone of positive semidefinite operators is also salient and full-dimensional.
\item
	The space of positive functions:

	Let $\cX$ be an arbitrary set, then the space of all functions $\cX \to \rset$ gives rise to a directed ordered vector space with the product order 
\[
	f \geq g \iff f(x) \geq g(x) \text{ for all } x \in \cX.
\]
	That is, the convex cone of componentwise positive functions is also salient and full-dimensional. 
	Let $\cX$ be a locally compact Hausdorff space, then this product order also applies to the spaces $C(\cX, \rset)$ of continuous functions and $C_c(\cX, \rset)$ of compactly supported continuous functions.

\exmsymbol
\end{enumerate}
\end{exms}

	If $V$ is a directed preordered space, then one can always find a basis of $V$ consisting solely of positive elements. 
	The argument is as follows.
	Let $\{b_i\}_{i \in I}$ be a Hamel basis of $V$, then for each $i \in I$, there exist $b_i^+, b_i^- \in V_+$ such that $b_i = b_i^+ - b_i^-$. 
	Now the collection $\{b_i^+, b_i^- \mid i \in I\}$ spans the vector space $V$ and can be reduced to a basis of $V$.
	
\begin{definition}
	Let $V$ be a preordered vector space and $E \subseteq V$ a linear subspace. 
	$E$ is said to be \emph{cofinal}, if for every $v \in V$ there exists an $x \in E$ with $x \geq v$.
\end{definition}

	There are multiple useful characterizations for a subspace to be cofinal. 
	We list them in the following lemma which is related to \cite[Prop.\ 34.1]{choquet69}, and \cite[Lem.\ 2.1]{didio24TsystemsArxiv}.
\begin{lem}
	Let $V$ be a preordered vector space and $E \subseteq V$ a linear subspace. 
	Then the following are equivalent:
\begin{enumerate}[(i)]
\item
	$E$ is cofinal.
\item
	For every $v \in V$ there exists an $x \in E$ with $x \geq v$.
\item
	For every $v \in V$ there exists an $x \in E$ with $v \geq x$.
\item
	$V = E + V_+$.
\end{enumerate}
\end{lem}
\begin{proof}
(i) $\iff$ (ii):
	This is simply the definition of cofinality.

(ii) $\Rightarrow$ (iii): 
	Let $v \in V$, then there exists an $x \in E$ such that $-x \geq -v$. 
	Hence, 
\[
	-x - (-v) = v - x \geq 0
\]
	and $v \geq x \in E$.

(iii) $\Rightarrow$ (iv):
	Let $v \in V$, then there exists an $x \in E$ such that $v \geq x$.
	Hence, $v - x \in V_+$ and $v = x + (v-x) \in E + V_+$. 
	Therefore $V \subseteq E + V_+$ and because $E, V_+ \subseteq V$, it follows, that $V = E + V_+$.

(iv) $\Rightarrow$ (ii):
	Let $v \in V$, then there exist $x \in E$ and $c \in V_+$ with $-v = -x + c$. 
	Hence, $x - v = c \geq 0$ and $x \geq v$.
\end{proof}

\begin{rem}
	Let $V, W$ be real preordered vector spaces, then $V \otimes W$ is a real preordered vector space with positivity cone 
\[
	V_+ \otimes W_+ := \conv\{v \otimes w \mid v \in V_+, w \in W_+\}.
\]
	$V_+ \otimes W_+$ is by definition a convex cone and is called the minimal tensor product of cones. 
	If $V_+ \subseteq V$ and $W_+ \subseteq W$ are full-dimensional, then $V_+ \otimes W_+ \subseteq V \otimes W$ is also full-dimensional. 
	This is simply because 
\[
	v \otimes w = v^+ \otimes w^+ + v^- \otimes w^- - v^+ \otimes w^- - v^- \otimes w^+
\]
	for all $v \in V$, $w \in W$ where $v = v^+ - v^-$ and $w = w^+ - w^-$ with $v^+, v^- \in V_+$, $w^+, w^- \in W_+$. 
	Hence, this defines a tensor product of preordered vector spaces such that the tensor product of two directed preordered vector spaces is again a directed preordered vector space. 
	Additionally, if $V_+ \subseteq V$ is full-dimensional and $E \subseteq W$ is cofinal, then $V \otimes E \subseteq V \otimes W$ is also cofinal, see the following lemma.
\exmsymbol
\end{rem}

\begin{lem}[see \protect{\cite[Lem.\ 5.4]{langer26_phd} or \cite[Lem.\ 2.10]{didiolang26_matrixpospres}}]
\label{lem:FullDimensionalConeExtendsCofinality}
	Let $V, W$ be real preordered vector spaces. 
	Let $V$ be a directed set and let $E \subseteq W$ be a cofinal subspace. 
	Then $V \otimes E \subseteq V \otimes W$ is a cofinal subspace with respect to the positivity cone $V_+ \otimes W_+$.
\end{lem}

	Positivity cones are used to define positive functionals.

\begin{definition}
	Let $V$ be a preordered vector space and $L: V \to \rset$ be a linear functional. 
	$L$ is said to be \emph{positive}, if $L(v) \geq 0$ for all $v \in V_+$.
\end{definition}

	Every positive functional on a cofinal subspace can be extended to a positive functional on the whole space. 
	This is the M.\ Riesz extension theorem, given as follows.

\begin{thm}[M.\ Riesz extension theorem, see \protect{\cite[34.2]{choquet69}, \cite[Thm.\ 2.3]{didio24TsystemsArxiv} or \cite[Prop.\ 1.7]{schmudMomentBook}}]
\label{thm:MRieszExtensionTheorem}
	Let $V$ be a real preordered vector space, let $E \subseteq V$ be a cofinal subspace, and let $L: E \to \rset$ be a positive linear functional.
	Then there exists a linear extension $\widetilde{L}: V \to \rset$ such that
\begin{enumerate}[(i)]
\item
	 $\widetilde{L}$ is a positive linear functional and
\item
	$\widetilde{L}|_E = L$.
\end{enumerate}
\end{thm}

	Note, that the M.\ Riesz extension theorem (\Cref{thm:MRieszExtensionTheorem}) needs a real vector space $V$. 
	Hence, why this paper will mostly work with real vector spaces.


\subsection{The Riesz--Markov--Kakutani Representation Theorem}
\label{subsec:RieszMarkovKakutani}

\begin{thm}[Riesz--Markov--Kakutani representation theorem, see \protect{\cite[Thm.\ 6.19]{rudin87}}]
\label{thm:RieszMarkovKakutani:FiniteMeasures}
	Let $\cX$ be a locally compact Hausdorff space, and let $L: C_c(\cX, \rset) \to \rset$ be a linear functional such that
\begin{enumerate}[(i)]
\item
	$L$ is positive, i.e., $L(f) \geq 0$ for all $f \in C_c(\cX, \rset)_+$ and 
\item
	$L$ is bounded, i.e., there exists a constant $C \in \rset$ with $L(f) \leq C \|f\|_\infty$ for all $f \in C_c(\cX, \rset)$.
\end{enumerate} 
	Then there exists a unique positive finite regular Borel measure $\mu: \Bor(\cX) \to \rset$ such that $\mu$ represents $L$, i.e., 
\[
	L(f) = \int_\cX f ~\diff\mu
\]
	for all $f \in C_c(\cX, \rset)$.
\end{thm}
	When $\cX$ is $\sigma$-compact, we also get the following integral approximation.

\begin{lem}
\label{lem:IntegralApproximationCc}
	Let $\cX$ be a $\sigma$-compact locally compact Hausdorff space and $\mu$ a positive Borel measure on $\cX$. 
	Then, for each positive continuous function $f \in C(\cX, \rset)_+$, 
\[
	\int_\cX f ~\diff\mu = \sup\left\{\int_\cX \phi ~\diff\mu ~\middle|~ \phi \in C_c(\cX, \rset)_+, \phi \leq f\right\}.
\]
\end{lem}
\begin{proof}
	By the monotonicity of the integral, 
\[
	\int_\cX \phi ~\diff \mu \leq \int_\cX f ~\diff\mu
\]
	for all measurable $0 \leq \phi \leq f$. 
	Hence 
\[
	\sup\left\{\int_\cX \phi ~\diff\mu ~\middle|~ \phi \in C_c(\cX, \rset)_+, \phi \leq f\right\} \leq \int_\cX f ~\diff\mu.
\]
	Since $\cX$ is a locally compact Hausdorff space, by Urysohn's lemma, there exists for each compact set $K \subseteq \cX$ a continuous compactly supported function $f \in C_c(\cX, \rset)$ such that $0 \leq f \leq 1$ and $f \equiv 1$ on $K$, see e.g. \cite[Lem.\ 2.12]{rudin87}. 
	Since $\cX$ is $\sigma$-compact, there exists an increasing sequence $K_1 \subseteq K_2 \subseteq ...$ of compact sets such that 
\[
	\cX = \bigcup_{i = 1}^\infty K_i = \lim_{i \to \infty} K_i.
\]
	Hence, for every $i \in \nset$, there exists a function $g_i \in C_c(\cX, \rset)$ with $0 \leq g_i \leq 1$ and $g_i \equiv 1$ on $K_i$. 
	Define $f_i := g_i \cdot f$, then $f_i \in C_c(\cX, \rset)_+$, $f_i \leq f$, and $f_i \overset{i \to \infty}{\longrightarrow} f$ pointwise and monotonous from below. 
	By the monotone convergence theorem, 
\[
	\int_\cX f ~\diff\mu = \lim_{i \to \infty} \int_\cX f_i ~\diff\mu
\]
	and thus 
\[
	\int_\cX f ~\diff\mu = \sup\left\{\int_\cX \phi ~\diff\mu ~\middle|~ \phi \in C_c(\cX, \rset)_+, \phi \leq f\right\}.
\]
\end{proof}

	Note, that without $\sigma$-compactness \Cref{lem:IntegralApproximationCc} can fail. 
	See the following example.
\begin{exm}
	Let $\cX := \rset$ be equipped with the discrete topology, i.e., every set is open in $\cX$. 
	This space $\cX$ is locally compact, Hausdorff but not $\sigma$-compact.
	A set is compact in this topology, if and only if, it is finite. 
	Hence, define the positive Borel measure $\mu: \Bor(\cX) \to \rset$ with 
\[
	\mu(S) := \begin{cases}
		0, & S\text{ is finite or countable}
	\\	\infty, & \text{else}
	\end{cases}
\]
	for $S \subseteq \rset$. 
	The $\sigma$-additivity of $\mu$ follows from countable unions of countable sets being countable again. 
	Hence for every $\phi \in C_c(\cX, \rset)$, 
\[
	\int_\cX \phi ~\diff\mu = 0.
\]
	Now let $f \equiv 1$ be the indicator function of the whole set $\cX$. 
	Then 
\[
	\int_\cX f ~\diff\mu = \mu(\cX) = \infty \neq 0 = \sup\left\{\int_\cX \phi ~\diff\mu ~\middle|~ \phi \in C_c(\cX, \rset)_+, \phi \leq f\right\}.
\]
\exmsymbol
\end{exm}


\subsection{Vector-Valued Integration and Pettis' Measurability Theorem}
\label{subsec:VectorValuedIntegrationPettisTheorem}

	Let $V$ be a real or complex vector space. 
	For real valued functions measurability is defined clearly, but there are multiple ways of defining a notion of measurability for vector-valued functions $f: \cX \to V$ on some measurable space $\cX$. 
	Pettis used functionals from the dual space $V'$ to define weak measurability while, in the case that $V$ is a Banach space, Bochner used the Banach topology to define strong measurability. 
	Concretely, they are defined as follows.
\begin{definition}
	Let $(\cX, \cA, \mu)$ be a measure space and $B$ a real or complex Banach space. 
	A function $f: \cX \to B$ is called
\begin{itemize}
\item
	a \emph{simple function}, if there exists an $n \in \nset$, vectors $x_1, \dots, x_n \in B$ and sets $\Omega_1, \dots, \Omega_n \in \cA$ such that 
\[
	f = \sum_{k = 1}^n x_k \one_{\Omega_k}.
\]
\item
	\emph{countably valued}, if $\im(f)$ is a countable set.
\item
	\emph{Bochner measurable} or \emph{strongly measurable}, if there exists a $\mu$-null set $N \in \cA$ and a sequence of simple functions $(s_i)_{i \in \nset}: \cX \to B$ such that 
\[
	\lim_{i \to \infty} \|f(x) - s_i(x)\|_B = 0
\]
	for all $x \in \cX \setminus N$.
\item
	\emph{weakly measurable}, if $g \circ f: \cX \to \kset$ is measurable for all continuous linear functionals $g: B \to \kset$ where $\kset$ is the underlying field of the Banach space $B$.
\item
	\emph{$\mu$-almost surely separably valued}, if there exists a $\mu$-null set $N \in \cA$ such that $f(\cX \setminus N)$ is a separable topological space with the subspace topology induced by $B$.
\end{itemize}
\end{definition}

	Pettis then characterized the connection between weak measurability and Bochner measurability.

\begin{thm}[Pettis' measurability theorem, see \protect{\cite[Thm.\ 1.1]{pettis38}}]
\label{thm:PettisMeasurabilityTheorem}
	Let $(\cX, \cA, \mu)$ be a measure space, let $B$ be a Banach space, and let $f: \cX \to B$ be a function. 
	Then the following are equivalent:
\begin{enumerate}[(i)]
\item
	$f$ is Bochner measurable.
\item
	$f$ is weakly measurable and $\mu$-almost surely separably valued.
\end{enumerate}
\end{thm}

	The proof of Pettis' measurability theorem (\Cref{thm:PettisMeasurabilityTheorem}) can be altered such that the approximating sequence of simple functions for the Bochner measurable function $f$ has values only in $\im(f)$, the image of $f$. 
	This is the following lemma and the proof is almost identical to the proof of Pettis' measurability theorem, see \cite[Thm.\ 1.1]{pettis38}, \cite[Thm.\ E.9]{cohn13}, or \cite[Thm.\ 1.1.1]{arendtetal11}.
\begin{lem}
\label{lem:BochnerMeasurableApproximationInImage}
	Let $(\cX, \cA, \mu)$ be a measure space and $B$ be a Banach space.
	Let $f: \cX \to B$ be weakly measurable and $\mu$-almost surely separably valued, then there exists a $\mu$-null set $N \in \cA$ such that $f$ can be approximated pointwise on $\cX \setminus N$ by simple functions with values in $\im(f) \cup \{0\}$.
	Hence, $f$ is Bochner measurable.
\end{lem}
\begin{proof}
	Let $\Omega_0 \subseteq \cX$ be a $\mu$-null set such that $f(\cX \setminus \Omega_0)$ is separable. 
	Let $(y_n)_{n \in \nset}$ be a dense subset of $f(\cX \setminus \Omega_0)$ and define 
\[
	Z := \overline{\tspan\{y_n \mid n \in \nset\}}.
\]
	$Z \subseteq B$ is a separable Banach space with $f(\cX \setminus \Omega_0) \subseteq Z$. 
	Hence, there exist a dense subset $(z_k)_{k \in \nset}$ of $Z$ and by the Hahn--Banach theorem, there exist continuous linear functionals $\phi_k: B \to \rset$, $k \in \nset$ with 
\[
	\|\phi_k\|_{B'} = 1 \quad\text{and}\quad \phi_k(z_k) = \|z_k\|_B
\]
	for all $k \in \nset$. 
	Since $(z_k)_{k \in \nset} \subseteq Z$ is dense, 
\[
	\|z\|_B = \sup\{|\phi_k(z)| \mid k \in \nset\}
\]
	holds for all $z \in Z$, see also \cite[Lem.\ E.10]{cohn13}. 
	Since $f$ is weakly measurable, it follows that $x \mapsto |\phi_k(f(x)) - z|$ is measurable for all $k \in \nset$ and thus ${x \mapsto \|f(x) - z\|_B}$ is measurable.
	Now, since $(y_n)_{n \in \nset}$ is dense in $f(\cX \setminus \Omega_0)$, there exist for every $x \in \cX \setminus \Omega_0$ and every $\varepsilon > 0$, a $k \in \nset$ such that $\|f(x) - y_k\| < \varepsilon$.
	Define the sets 
\[
	\Omega_{k, \varepsilon} := \{x \in \cX \setminus \Omega_0 \mid \|f(x) - y_k\|_B < \varepsilon\}
\]
	and the functions $f_\varepsilon: \cX \to B$, 
\[
	f_\varepsilon := \sum_{k = 1}^\infty y_k \one_{\Omega_{k, \varepsilon}}.
\]
	Now $\im(f_\varepsilon) \subseteq \im(f) \cup \{0\}$ and $\|f-f_\varepsilon\|_B < \varepsilon$ for all $x \in \cX \setminus \Omega_0$.
	Additionally, because the maps $x \mapsto \|f(x) - z\|_B$ are measurable for all $z \in Z$, all $\Omega_{k, \varepsilon} \subseteq \cX$ for $k \in \nset$, $\varepsilon > 0$ are measurable. 
	Hence, $f_\varepsilon$ are measurable functions for $\varepsilon > 0$.
	Since $f_\varepsilon$ only has countably many values for $\varepsilon > 0$, it can be approximated by simple functions with the same values and thus $f$ can be approximated pointwise on $\cX \setminus \Omega_0$ by simple functions with values in $\im(f) \cup \{0\}$.
\end{proof}
	\Cref{lem:BochnerMeasurableApproximationInImage} also extends and simplifies the previous work of \cite[Lem.\ 4.14]{langer26_phd}.


\subsection{Adapted Spaces and Choquet's Basic Representation Theorem}
\label{subsec:AdaptedSpacesChoquetBasicRepresentationTheorem}

\begin{definition}[see e.g. \protect{%
	\cite{choquet69}, 
	\cite[Def.\ 2.7]{didio24TsystemsArxiv}, 
	\cite[Def.\ 2.4]{didio26_arxiv_linearoperators}, or 
	\cite[Def.\ 1.5]{schmudMomentBook}%
}]
	Let $\cX$ be a locally compact Hausdorff space, $E \subseteq C(\cX, \rset)$ be a vector space, and denote by $E_+ := E \cap C(\cX, \rset)_+$ the subspace of positive functions of $E$. 
	$E$ is called an \emph{adapted space}, if 
\begin{enumerate}[(i)]
\item
	$(E, E_+)$ is a directed preordered vector space, 
\item
	for all $x \in \cX$, there exists an $f \in E_+$ with $f(x) > 0$, 
\item
	for every $g \in E_+$, there exists an $f \in E_+$ which dominates $g$, i.e., such that for every $\varepsilon > 0$, there exists an $h_\varepsilon \in C_c(\cX, \rset)$ with $g \leq \varepsilon f + h_\varepsilon$.
\end{enumerate}
\end{definition}

\begin{exms}
\label{exms:AdaptedSpaces}
	In the following, we give some common examples of adapted spaces.
\begin{enumerate}[(a)]
\item
	The space of polynomials:

	Let $n \in \nset$ and let $\rset[x_1, \dots, x_n]$ be the space of polynomials in $n$ variables. 
	Then $\rset[x_1, \dots, x_n]$ is an adapted space on $\rset^n$. 
	Similarly, for every closed set $K \subseteq \rset^n$, $\rset[x_1, \dots, x_n]|_K$ is also an adapted space on $K$.
\item
	The space of exponential functions:

	Let $n \in \nset$ and let 
\[
	E = \tspan\{e^{\alpha_1 x_1 + \dots + \alpha_n x_n} \mid \alpha \in \rset^n\}
\]
	be the space of exponential functions on $\rset^n$. 
	Then $E \subseteq C(\rset^n, \rset)$ is an adapted space. 
	Similarly, for every sequence $0 = \alpha_0 < \alpha_1 < \alpha_2 < ...$, the space 
\[
	\tspan\{e^{\alpha_i x}, e^{-\alpha_i x} \mid i \in \nset_0\}
\]
	is an adapted space on $\rset$.
\end{enumerate}
\end{exms}

\begin{lem}[see e.g. \protect{%
	\cite{choquet69}, 
	\cite[Lem.\ 2.8]{didio24TsystemsArxiv},
	\cite[Lem.\ 2.5]{didio26_arxiv_linearoperators}, or 
	\cite[Lem.\ 1.6]{schmudMomentBook}%
}]
\label{lem:AdaptedSpaceCcDominated}
	Let $\cX$ be a locally compact Hausdorff space and $E \subseteq C(\cX, \rset)$ be an adapted space. 
	Then, for every $g \in C_c(\cX, \rset)$, there exists an $f \in E_+$ with $f > g$.
\end{lem}

\begin{thm}[Choquet's basic representation theorem, see \protect{%
	\cite[Thm.\ 34.6]{choquet69}, 
	\cite[Thm.\ 2.9]{didio24TsystemsArxiv}, 
	\cite[Thm.\ 2.7]{didio26_arxiv_linearoperators}, or 
	\cite[Thm.\ 1.8]{schmudMomentBook}%
}]
\label{thm:ChoquetBasicRepresentationTheorem}
	Let $\cX$ be a locally compact Hausdorff space, $E \subseteq C(\cX, \rset)$ be an adapted space, and $L: E \to \rset$ be a linear functional.
	Then the following are equivalent:
\begin{enumerate}[(i)]
\item 
	$L$ is a positive functional.
\item
	There exists a positive Borel measure $\mu: \Bor(\cX) \to \rset$ such that 
\[
	L(f) = \int_\cX f ~\diff\mu
\]
	for all $f \in E$.
\end{enumerate}
\end{thm}

	Note, that \Cref{thm:ChoquetBasicRepresentationTheorem} does not give a unique Borel measure. 

	From \Cref{thm:ChoquetBasicRepresentationTheorem} follows Haviland's theorem which is important in the theory of moments. 
	It characterizes, when a linear functional on the polynomials $L: \rset[x_1, \dots, x_n] \to \rset$ is represented by a moment measure $\mu$.

\begin{cor}[Haviland's Theorem, see \protect{\cite{havila36}}]
\label{cor:HavilandTheorem:Polynomials}
	Let $n \in \nset$, $K \subseteq \rset^n$ closed, and $L: \rset[x_1, \dots, x_n] \to \rset$ be a linear functional. 
	Denote by $\Pos_\rset(K)$ the set of polynomials $p \in \rset[x_1, \dots, x_n]$ which are positive on $K$, i.e., $p(x) \geq 0$ for all $x \in K$.
	Then the following are equivalent: 
\begin{enumerate}[(i)]
\item $L$ is positive on $K$, i.e., $L(p) \geq 0$ for all $p \in \Pos_\rset(K)$.
\item $(L(x^\alpha))_{\alpha \in \nset_0^n}$ is a $K$-moment sequence, i.e., there exists a measure 
	$\mu: \Bor(K) \to \rset$ such that 
\[
	L(p) = \int_K p(x) ~\diff\mu(x)
\]
	for all $p \in \rset[x_1, \dots, x_n]$.
\end{enumerate}
\end{cor}
\begin{proof}
	This follows from \Cref{thm:ChoquetBasicRepresentationTheorem} with $\cX := K$ and $E := \rset[x_1, \dots, x_n]|_K$.
\end{proof}

	Just like in Choquet's basic representation theorem, the measure $\mu$ from Haviland's theorem (\Cref{cor:HavilandTheorem:Polynomials}) is not unique.


\section{Map-Valued Measures}
\label{sec:MapValuedMeasures}

\begin{definition}
	Let $(\cX, \cA)$ be a measurable space and $V$ a real directed preordered vector space. 
	A function $\mu: \cA \to \cL(V, \rset)$ is a \emph{positive map-valued measure} if 
\[
	\mu[v] := \mu(\,\cdot\,)(v): \cA \to \rset
\]
	is a positive finite measure for all $v \in V_+$. 
	If $\cX$ is a topological space, then a positive Borel map-valued measure $\mu: \Bor(\cX) \to \cL(V, \rset)$ is called \emph{regular} if $\mu[v]$ is a regular Borel-measure for all $v \in V_+$.
\end{definition}

	Our definition of map-valued measures coincides with the definition of non-negative operator-valued measures from \cite[Sec.\ 2]{cimpzalar13}. 
	Cimprič and Zalar restricted the vector space $V$ to the space $B(\cH)_h$ of bounded self-adjoint operators on a Hilbert space $\cH$ while we aim to generalize the theories to arbitrary vector spaces $V$. 
	The boundedness of the measures $\mu[v]$ is needed to make sense of linear combinations $\mu[v] - \mu[w] = \mu[v-w]$ which will be signed measures in general. 
	Additionally, defining map-valued measures on positive elements $v \in V_+$ is enough, as $V$ is a directed set, i.e., $V_+ \subseteq V$ is full-dimensional. 
	Hence $\mu[v]$ linearly extends to all $v \in V$.

	We only define map-valued measures $\mu: \cA \to \cL(V, \rset)$ into the set of functionals on $V$.
	A generalization to arbitrary map-valued measures $\cA \to \cL(V, W)$ where $V, W$ are real vector spaces is possible too. 
	For the integration of $W$-valued functions, one would use the Pettis-integral, i.e., integration via functionals $W \to \rset$. 
	For this, the set of functionals on $W$ needs to separate points or else the integral will lack uniqueness.

	Integration with respect to map-valued measures is defined in the following way.

\begin{definition}
	Let $(\cX, \cA)$ be a measurable space, $V$ a real directed preordered vector space, and $\mu: \cA \to \cL(V, \rset)$ a positive map-valued measure. 
	Denote by $M(\cX, \rset)$ the set of measurable functions $\cX \to \rset$ and let $f \in V \otimes M(\cX, \rset)$, i.e., 
\[
	f = v_1 \otimes f_1 + \dots + v_n \otimes f_n
\]
	for some $n \in \nset$, $v_1, \dots, v_n \in V$ and $f_1, \dots, f_n \in M(\cX, \rset)$.
	Then $f$ is $\mu$-integrable, if $f_1, \dots, f_n$ is $\mu[v]$-integrable for all $v \in V_+$ and 
\[
	\int_\cX f ~\diff\mu := \int_\cX f_1 ~\diff\mu[v_1] + \dots + \int_\cX f_n ~\diff\mu[v_n].
\]
\end{definition}

	Integration with respect to map-valued measures does not depend on the representation of $p \in V \otimes M(\cX, \rset)$ as long as this representation is integrable, i.e., it does not depend on which $v_i$ and $f_i$ we choose.
\begin{lem}
	Let $(\cX, \cA)$ be a measurable space, $V$ a real directed preordered vector space, and $\mu: \cA \to \cL(V, \rset)$ a positive map-valued measure. 
	Denote by $M(\cX, \rset)$ the set of measurable functions $\cX \to \rset$ and let $f \in V \otimes M(\cX, \rset)$ with 
\[
	f 
	= v_1 \otimes f_1 + \dots + v_n \otimes f_k 
	= \widetilde{v}_1 \otimes \widetilde{f}_1 + \dots + \widetilde{v}_l \otimes \widetilde{f}_l
\]
	for some $k, l \in \nset$, $v_1, \dots, v_k, \widetilde{v}_1, \dots, \widetilde{v}_l \in V$, 
	and $f_1, \dots, f_k, \widetilde{f}_1, \dots, \widetilde{f}_l \in M(\cX, \rset)$. 
	Additionally, let $f$ be $\mu$-integrable with respect to both representations, then 
\[
	\sum_{i = 1}^k \int_\cX f_i(x) ~\diff\mu[v_i](x) = \sum_{i = 1}^l \int_\cX \widetilde{f}_i(x) ~\diff\mu[\widetilde{v}_i](x).
\]
\end{lem}
	The proof is adapted from \cite[Lem.\ 4.11]{langer26_phd}.
\begin{proof}
	Let $b_1, \dots, b_m \in V_+$, $m \in \nset$ be a basis of the finite-dimensional subspace 
\[
	\tspan\{v_1, \dots, v_k, \widetilde{v}_1, \dots, \widetilde{v}_l\} \subseteq V,
\]
	then there exist $\lambda_{i, j} \in \rset$ for $i = 1, \dots, k$, $j = 1, \dots, m$ such that 
\[
	v_i = \sum_{j = 1}^m \lambda_{i, j} b_j
\]
	for $i = 1, \dots, k$. 
	Hence, 
\[
	f 
	= v_1 \otimes f_1 + \dots + v_k \otimes f_k
	= b_1 \otimes g_1 + \dots + b_m \otimes g_m
\]
	with $g_j := \sum_{i = 1}^k \lambda_{i, j} f_i$ for $j = 1, \dots, m$. 
	By the linearity property of $\mu$, 
\begin{align*}
	\sum_{i = 1}^k \int_\cX f_i(x) ~\diff\mu[v_i](x) 
	&= \sum_{i = 1}^k\sum_{j = 1}^m \lambda_{i, j} \int_\cX f_i(x) ~\diff\mu[b_j](x) 
\\	&= \sum_{j = 1}^m \int_\cX g_j(x) ~\diff\mu[b_j](x)
\end{align*}
	where the integrability is ensured by $f_1, \dots, f_k$ being $\mu[v]$-integrable for all $v \in V_+$ since $f$ is $\mu$-integrable.
	Similarly, 
\[
	\sum_{i = 1}^l \int_\cX \widetilde{f}_i(x) ~\diff\mu[\widetilde{v}_i](x)
	= \sum_{j = 1}^m \int_\cX g_j(x) ~\diff\mu[b_j](x)
\]
	where we used the fact that $b_1, \dots, b_m$ is a basis to obtain the polynomials $g_1, \dots, g_m$ in the integral again.
\end{proof}

	Integration with respect to map-valued measures also preserves positivity, see the following theorem.

\begin{thm}[positivity of map-valued integrals]
\label{thm:PositivityMapValueIntegrals}
	Let $(\cX, \cA)$ be a measurable space, $V$ be a real preordered vector space, and $\mu: \cA \to \cL(V, \rset)$ a positive map-valued measure. 
	Then for all $p \in V \otimes M_s$ such that $\im(p) \subseteq V_+$, 
\[
	\int_\cX p ~\diff\mu \geq 0
\]
	where 
\[
	M_s := \{g: \cX \to \rset \mid g\text{ measurable and }\mu\text{-almost surely separably valued}\}.
\]
	The integral here is allowed to diverge to $\infty$ which is understood as positive, i.e., $\infty \geq 0$.
\end{thm}
	This proof is adapted from \cite[Lem.\ 4.15]{langer26_phd}.
\begin{proof}
	Let $p \in V \otimes M_s$ such that $\im(p) \subseteq V_+$, then $p$ is a function $\cX \to V$.
	Define 
\[
	C := \conv\{\lambda x \mid \lambda \in \rset_{\geq 0}, x \in \im(p)\} \subseteq V_+, 
	\quad
	E := \tspan\{\im(p)\}
\]
 	then $E \subseteq V$ is a finite-dimensional vector space, $C \subseteq E$ is full-dimensional, and $\im(p) \subseteq E$. 
	Our goal is to approximate $p$ with simple functions $s_i: \cX \to E$ on some basis of $E$. 
	We equip $E$ with the Euclidean topology to make sense of this approximation. 
	Since $p \in V \otimes M_s$ is a finite sum, $p: \cX \to E$ is weakly measurable and $\mu$-almost surely separably valued.
	Pettis' measurability theorem (\Cref{thm:PettisMeasurabilityTheorem}) now asserts that $p$ can be approximated almost surely pointwise by simple functions $s_i: \cX \to E$, $i \in \nset$. 
	By \Cref{lem:BochnerMeasurableApproximationInImage}, these simple functions can be chosen to be positive, i.e., $\im(s_i) \subseteq \im(p) \cup \{0\} \subseteq V_+$ for all $i \in \nset$. 
	Let $b_1, \dots, b_n \in E$ be a basis of $E$. 
	Since $C \subseteq E$ is full-dimensional, $b_1, \dots, b_n$ can be chosen to be positive elements. 
	Now, decompose $p = b_1 \otimes p_1 + \dots + b_n \otimes p_n$ and $s_i = b_1 \otimes f_{i, 1} + \dots + b_n \otimes f_{i, n}$ over the basis $b_1, \dots, b_n$. 
	The functions $f_{i, 1}, \dots, f_{i, n}: \cX \to \rset$, $i \in \nset$ are measurable simple functions.
	Hence, by the monotone convergence theorem, 
\begin{align*}
	\lim_{i \to \infty} \left|\int_\cX p - s_i ~\diff\mu\right|
	&= \lim_{i \to \infty} \left|\sum_{k = 1}^n \int_\cX p_k - f_{k, i} ~\diff\mu[b_k]\right|
\\	&\leq \lim_{i \to \infty} \sum_{k = 1}^n \int_\cX |p_k - f_{k, i}| ~\diff\mu[b_k]
\\	&= 0
\end{align*}
	and 
\[
	\lim_{i \to \infty} \int_\cX s_i ~\diff\mu = \int_\cX p ~\diff\mu.
\]
	Since $s_i$ are positive simple functions, the integrals are positive as well, i.e., 
\[
	\int_\cX s_i ~\diff\mu \geq 0
\]
	for all $i \in \nset$ and thus 
\[
	\int_\cX p ~\diff\mu \geq 0.
\]
\end{proof}

	The Riesz--Markov--Kakutani theorem (\Cref{thm:RieszMarkovKakutani:FiniteMeasures}) also has an extension to map-valued measures, see the following theorem.

\begin{thm}[Riesz--Markov--Kakutani theorem for tensor spaces]
\label{thm:RieszMarkovKakutaniTensorSpaces}
	Let $\cX$ be a locally compact Hausdorff space, $V$ a real directed preordered vector space, and let $L: V \otimes C_c(\cX, \rset) \to \rset$ be a linear functional such that
\begin{itemize}
\item
	$L$ is positive, i.e., $L(f) \geq 0$ for all $f \in V_+ \otimes C_c(\cX, \rset)_+$ and
\item
	For all $v \in V_+$, there exists a constant $C_v \geq 0$ such that 
\[
	L(v \otimes f) \leq C_v \|f\|_\infty
\]
	for all $f \in C_c(\cX, \rset)$.
\end{itemize}
	Then there exists a unique positive regular Borel map-valued measure $\mu: \Bor(\cX) \to \cL(V, \rset)$ such that 
\[
	L(f) = \int_\cX f ~\diff\mu
\]
	for all $f \in V \otimes C_c(\cX, \rset)$.
\end{thm}
	This proof is an adapted version of \cite[Prop.\ 2]{cimpzalar13}.
\begin{proof}
	For $v \in V$, define the linear functionals $L_v: C_c(\cX, \rset) \to \rset$ 
\[
	L_v(f) := L(v \otimes f)
\]
	for $f \in C_c(\cX, \rset)$. 
	Then $L_v$ is a positive bounded functional for $v \in V_+$.
	Hence, by the Riesz--Markov--Kakutani representation theorem (\Cref{thm:RieszMarkovKakutani:FiniteMeasures}), there exist unique positive finite Borel measures $\mu_v: \Bor(\cX) \to \rset$, such that 
\[
	L_v(f) = L(v \otimes f) = \int_\cX f(x) ~\diff\mu_v(x)
\]
	for all $f \in C_c(\cX, \rset)$. 
	The uniqueness of the measures from the Riesz--Markov--Kakutani representation theorem (\Cref{thm:RieszMarkovKakutani:FiniteMeasures}) implies, 
\[
	\mu_{v+w} = \mu_v + \mu_w
\]
	for all $v, w \in V_+$.
	Let $v^+, v^-, w^+, w^- \in V_+$ with $v^+ - v^- = w^+ - w^-$, then 
\[
	\mu_{v^+} + \mu_{w^-} = \mu_{w^+} + \mu_{v^-}
\]
	and 
\[
	\mu_{v^+} - \mu_{v^-} = \mu_{w^+} - \mu_{w^-}.
\]
	This ensures the well-definedness of the definition 
\[
	\mu_v := \mu_{v^+} - \mu_{v^-}
\]
	for $v \in V$ where $v^+, v^- \in V_+$ with $v = v^+ - v^-$.
	Hence, the map-valued measure $\mu: \Bor(\cX) \to \cL(V, \rset)$, $\mu[v] := \mu_v$ is a positive regular Borel measure such that 
\[
	L(p) = \int_\cX p ~\diff\mu
\]
	for all $p \in V \otimes C_c(\cX, \rset)$.
\end{proof}


\section{Haviland's Theorem on Tensor Spaces}
\label{sec:HavilandTheoremTensorSpaces}

\begin{thm}[generalized Haviland theorem]
\label{thm:Haviland:TensorCones}
	Let $V$ be a real directed preordered vector space, $\cX$ be a locally compact $\sigma$-compact separable Hausdorff space, and $E \subseteq C(\cX, \rset)$ be an adapted subspace such that $1 \in E$.
	Define 
\[
	\Pos := \{p \in V \otimes E \mid p(x) \in V_+ \text{ for all }x \in \cX\}
\]
	For a linear functional $L: V \otimes E \to \rset$ the following are equivalent:
\begin{enumerate}[(i)]
\item
	$L$ is positive on $\Pos$, i.e., $L(p) \geq 0$ for all $p \in \Pos$.
\item
	There exists a positive Borel map-valued measure $\mu: \Bor(\cX) \to \cL(V, \rset)$ such that 
\[
	L(p) = \int_\cX p ~\diff\mu
\]
	for all $p \in V \otimes E$.
\end{enumerate}
\end{thm}
	The proof is adapted from the proof of Cimprič and Zalar's version of Haviland's theorem, see \cite[Thm.\ 3]{cimpzalar13}, as well as the proof of Choquet's basic representation theorem (\Cref{thm:ChoquetBasicRepresentationTheorem}).
\begin{proof}
	(ii) $\Rightarrow$ (i): 
	This is \Cref{thm:PositivityMapValueIntegrals}. 
	Note that all functions in $E \subseteq C(\cX, \rset)$ are continuous, hence all $p \in E$ are measurable with respect to the Borel $\sigma$-algebra on $\cX$.
	Since $\cX$ is separable and the image of a separable set under continuous maps is again separable, every $p \in E$ is measurable and separably valued. 

	(i) $\Rightarrow$ (ii):
	Define
\[
	F := \{f \in C(\cX, \rset) \mid |f| \leq p\text{ for some }p \in E_+\}.
\]
	Then $F$ is a real directed preordered vector space with positivity cone $F_+ := F \cap C(\cX, \rset)_+$ and $E \subseteq F$ is cofinal. 
	Since $V$ is also a real directed preordered vector space, by \Cref{lem:FullDimensionalConeExtendsCofinality}, $V \otimes E \subseteq V \otimes F$ is cofinal with respect to the full-dimensional positivity cone $V_+ \otimes F_+$. 
	Since, $L$ is positive on $\Pos$ and 
\[
	(V \otimes E) \cap (V_+ \otimes F_+) \subseteq \Pos,
\]
	$L$ is positive on $V \otimes E$.
	Hence, by the M.\ Riesz extension theorem (\Cref{thm:MRieszExtensionTheorem}), $L$ can be extended to a positive linear functional $\widetilde{L}: V \otimes F \to \rset$. 
	By \Cref{lem:AdaptedSpaceCcDominated}, $C_c(\cX, \rset) \subseteq F$ and for all $v \in V_+$ and $f \in C_c(\cX, \rset)$, 
\[
	v \otimes f \leq v \otimes |f| \leq \|f\|_\infty \cdot (v \otimes 1) \in V \otimes E.
\]
	Hence, for all $v \in V_+$, and $f \in C_c(\cX, \rset)$.
\[
	\widetilde{L}(v \otimes f) \leq L(v \otimes 1) \cdot \|f\|_\infty.
\]
	From the Riesz--Markov--Kakutani theorem (\Cref{thm:RieszMarkovKakutaniTensorSpaces}) now follows, that there exists a positive regular Borel map-valued measure $\mu: \Bor(\cX) \to \cL(V, \rset)$ such that 
\[
	\widetilde{L}(p) = \int_\cX p ~\diff\mu
\]
	for all $p \in V \otimes C_c(\cX, \rset)$.
	Now let $v \in V_+$ and $f \in E_+$, then $\mu[v]$ is a positive Borel measure.
	Since $\cX$ is $\sigma$-compact, by \Cref{lem:IntegralApproximationCc}, 
\begin{align*}
 	\int_\cX f(x) ~\diff\mu[v](x)
	&= \sup\left\{\int_\cX \phi(x) ~\diff\mu[v](x) \mid \phi \in C_c(\cX, \rset), \phi \leq f\right\}
\\	&= \sup\left\{\widetilde{L}(v \otimes \phi) \mid \phi \in C_c(\cX, \rset), \phi \leq f\right\}
\\	&\leq \widetilde{L}(v \otimes f) 
\\	&= L(v \otimes f).
\end{align*}
	Hence $f$ is $\mu[v]$-integrable for all $v \in V_+$ and thus every $p \in V \otimes E$ is $\mu$-integrable.
	Now let $v \in V_+$ and define the linear functionals $K_v: F \to \rset$, 
\[
	K_v(f) := \widetilde{L}(v \otimes f) - \int_\cX f(x) ~\diff\mu[v](x)
\]
	for all $f \in E$. 
	Then $K_v$ is positive on $E_+$ and $K_v(f) = 0$ for $f \in C_c(\cX, \rset)$. 
	Since $E$ is an adapted space, for every $g \in E_+$, there exists an $f \in E_+$ such that for every $\varepsilon > 0$, there exist a $h_\varepsilon \in C_c(\cX, \rset)$ with $g \leq \varepsilon f + h_\varepsilon$. 
	Hence, 
\[
	0 \leq K_v(g) \leq \varepsilon K_v(f) \overset{\varepsilon \to 0}{\longrightarrow} 0, 
\]
	i.e., $K_v = 0$ on $E$ for all $v \in V_+$ and 
\[
	\widetilde{L}(p) = L(p) = \int_\cX p ~\diff\mu
\]
	for all $p \in V \otimes E$.
\end{proof}

Cimprič and Zalar's extension of Haviland's theorem follows from \Cref{thm:Haviland:TensorCones}.

\begin{cor}[Cimprič and Zalar's version of Haviland's theorem, see \protect{\cite[Thm.\ 3]{cimpzalar13}}]
\label{cor:HavilandTheorem:OperatorPolynomials}
	Let $\cH$ be a Hilbert space and denote by $B(\cH)_h$ the space of bounded self-adjoint linear operators on $\cH$.
	Let $n \in \nset$, $K \subseteq \rset^n$ be a closed set, and define 
\[
	\Pos := \{p \in B(\cH)_h \otimes \rset[x_1, \dots, x_n] \mid p(x) \succeq 0\text{ for }x \in K\}
\]
	where $\succeq$ denotes positive semidefiniteness, i.e., the Loewner order.
	For a linear functional $L: B(\cH)_h \otimes \rset[x_1, \dots, x_n] \to \rset$, the following statements are equivalent:
\begin{enumerate}[(i)]
\item
	$L$ is positive on $\Pos$, i.e., $L(p) \geq 0$ for all $p \in \Pos$.
\item
	There exists a positive Borel map-valued measure $\mu: \Bor(\cX) \to \cL(B(\cH)_h, \rset)$ such that 
\[
	L(p) = \int_\cX p ~\diff\mu
\]
	for all $p \in B(\cH)_h \otimes \rset[x_1, \dots, x_n]$.
\end{enumerate}
\end{cor}
\begin{proof}
	Let $V := B(\cH)_h$, then $V$ is a real Banach space and a directed ordered vector space, see \Cref{exms:OrderedVectorSpaces}.(a). 
	Furthermore, $E := \rset[x_1, \dots, x_n]|_K$ is an adapted space which includes $1 = x^0$, see \Cref{exms:AdaptedSpaces}.(a). 
	Since $V_+ \otimes E_+ \subseteq \Pos$, every functional which is positive on $\Pos$ is also positive on $V_+ \otimes E_+$.
	Hence, \Cref{thm:Haviland:TensorCones} and \Cref{thm:PositivityMapValueIntegrals} complete the proof.
\end{proof}

	The space $\cX$ in \Cref{thm:Haviland:TensorCones} is required to be locally compact, $\sigma$-compact, separable, and Hausdorff. 
	This still includes most common spaces in mathematics, as every locally compact second-countable Hausdorff-space fulfills these requirements which includes every separable metrizable space.

\begin{exm}
	Let $V = \rset^{m \times m}$ be the space of $m \times m$ real matrices. 
	Equip $V$ with the product order or componentwise order, 
\[
	A \geq B \iff a_{i, j} \geq b_{i, j} \text{ for all }i, j = 1, \dots, m
\]
	where $A = (a_{i, j})_{i, j = 1}^m, B = (B_{i, j})_{i, j = 1}^m \in \rset^{m \times m}$. 
	This turns $V$ into a directed ordered vector space.
	It is easy to see, that a linear functional $L: V \otimes \rset[x_1, \dots, x_n] \to \rset$ is positive, if and only if, the linear functionals $L_{i, j}: \rset[x_1, \dots, x_n] \to \rset$ with 
\[
	L_{i, j}(p) = L(p \otimes E_{i, j})
\]
	are positive where $E_{i, j}$ is the $m \times m$ matrix with $1$ in the $(i, j)$-th component and $0$ elsewhere. 
	Hence, one could use Haviland's theorem on every component and get a matrix-valued measure to write $L$ as an integral. 
	This establishes a connection between positive matrix-valued measures and positive linear functionals on $V \otimes \rset[x_1, \dots, x_n]$. 
	This same connection is also given by \Cref{thm:Haviland:TensorCones}. 
	Hence, \Cref{thm:Haviland:TensorCones} generalizes both Cimprič and Zalar's Haviland theorem and this rather simple connection for matrix polynomials. 
\end{exm}


\section*{AI Use Disclosure}

This paper is written entirely by humans. 
Nonetheless, AI was used to peer review this version and to spot errors.


\section*{Funding}

The author is supported by the Deutsche Forschungs\-gemein\-schaft DFG with the grant DI-2780/2-1.


\bibliographystyle{amsalpha}

\providecommand{\bysame}{\leavevmode\hbox to3em{\hrulefill}\thinspace}
\providecommand{\MR}{\relax\ifhmode\unskip\space\fi MR }
\providecommand{\MRhref}[2]{%
  \href{http://www.ams.org/mathscinet-getitem?mr=#1}{#2}
}
\providecommand{\href}[2]{#2}

\end{document}